\documentclass[12pt,leqno]{article}
\usepackage{amsthm, amsmath, amssymb, amsfonts}

\input epsf.tex
\newdimen\epsfxsize
\newdimen\epsfysize
\def\qed{\vrule height5pt width3pt depth.5pt}

\theoremstyle{plain}
\newtheorem{thm}{Theorem}[section]
\newtheorem{cor}[thm]{Corollary}
\newtheorem{lem}[thm]{Lemma}
\newtheorem{prop}[thm]{Proposition}

\newtheorem{rem}{Remark}[section]

\begin{document}

\title{Pure Virtual Braids \\ Homotopic to the Identity Braid}

% information for author

\author{H. A. Dye \\
McKendree University \\
701 College Road \\
Lebanon, IL 62254 \\
hadye@mckendree.edu }

\maketitle

\begin{abstract} Two virtual link diagrams are homotopic if one may be transformed into the other by a sequence of virtual Reidemeister moves, classical Reidemeister moves, and self crossing changes. We recall the pure virtual braid group. We then describe the set of pure virtual braids that are homotopic to the identity braid.
\end{abstract}

\section{Introduction}

A virtual link diagram is a decorated immersion of $ n $ copies of $ S^1 $ with two types of crossings: classical and virtual. Classical crossings are indicated by over/under markings and virtual crossings are indicated by a solid encircled X. 
\begin{figure}[htb] \epsfysize = 0.75 in
\centerline{\epsffile{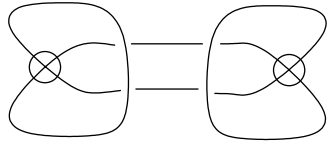}}
\caption{Kishino's knot}
\label{fig:kish}
\end{figure}
\noindent
An example of a virtual link diagram is shown in figure \ref{fig:kish}.

\begin{figure}[htb] \epsfysize = 0.5 in
\centerline{\epsffile{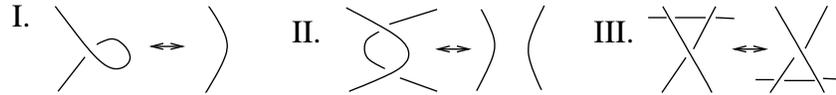}}
\caption{Classical Reidemeister moves}
\label{fig:rmoves}
\end{figure}

\begin{figure}[htb] \epsfysize = 1.25 in
\centerline{\epsffile{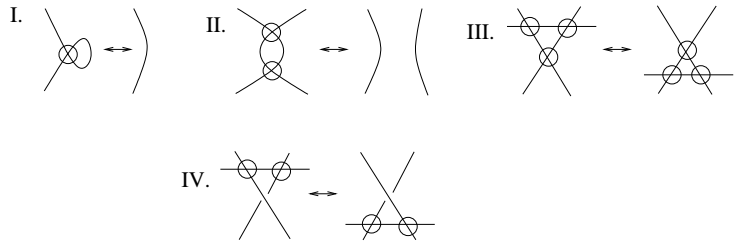}}
\caption{Virtual Reidemeister moves}
\label{fig:vrmoves}
\end{figure}

Virtual link theory is a generalization of classical knot theory that was introduced by Louis H. Kauffman in 1996 \cite{kvirt}.
Two virtual link diagrams are said to be \emph{equivalent} if one may be transformed into another by a sequence of classical Reidemeister moves (shown in figure \ref{fig:rmoves}) and virtual Reidemeister moves (shown in figure \ref{fig:vrmoves}). Classical link diagrams contain no virtual crossings and form a subset of the virtual link diagrams.

A \emph{representation} of a virtual link diagram is a pair $ (F,L) $ where $L$ is a link diagram on a closed, two dimensional surface $F$ taken up to Dehn twists and handle ($ S^1 \times I$)  additions and cancellations. Reidemeister moves can be performed on the surface. 
\begin{rem}
Detailed descriptions of representations are given in \cite{dk} or \cite{kvirt}. Abstract surfaces (punctured representations) are described in \cite{kamada2}.\end{rem}

We recall the following theorems:
\begin{thm}\label{corres} Classes of representations  are in one to one correspondence with equivalence classes of virtual link diagrams.\end{thm}
\textbf{Proof:} See \cite{kamada2} and \cite{kvirt}.\qed

\begin{thm}[Kuperberg]\label{kup} Representations of virtual link diagrams have a unique representative embedding class in the minimal genus surface that can support the diagram. \end{thm} 
\textbf{Proof:} See \cite{kup}.\qed

Recalling \cite{d-khom},
two virtual link diagrams are defined to be \emph{homotopic} if one diagram may be transformed into the other by a sequence of virtual Reidemeister moves, classical Reidemeister moves, and self crossing change. (By self crossing change, we mean changing the over/under markings at a crossing between two segments of the same link component.)

\begin{figure}[htb] \epsfysize = 1.25 in
\centerline{\epsffile{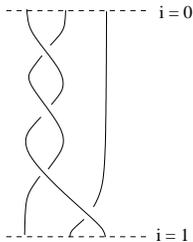}}
\caption{Example: braid diagram}
\label{fig:braidsample}
\end{figure}

In this paper, we focus on virtual braids (see \cite{kamadabraid} and \cite{sophia}) and pure virtual braids. A \emph{n-strand virtual braid diagram} is a decorated immersion of $ n $ copies of $ [0,1] $ into the plane. Let $ \lbrace f_1, f_2, \ldots f_n \rbrace $ denote the $n$ components. The set of endpoints $ \lbrace f_1 (i) , f_2 (i), \ldots
f_n (i) | i \in \lbrace 0,1 \rbrace \rbrace $ are contained on a line for each $i$. We refer to the set of points where $i=0$ as the upper endpoints and the other set ($i=1 $) as the lower endpoints as shown in figure \ref{fig:braidsample}.
Two virtual braid diagrams are said to \emph{virtually homotopic} if one can be transformed into the other by a sequence of Reidemeister moves, virtual Reidemeister moves, and self-crossing changes that leave the endpoints fixed.
 A virtual braid is \emph{pure} if the braid is equivalent to a braid with no self crossings, that is every crossing involves two distinct components. 
For a pure virtual braid, the order of the set of upper endpoints is the same as the order of the set of lower endpoints. 
We denote the set of n-strand pure virtual braids as 
$ VP_n $ and the set of n-strand pure classical braids as $BP_n$. We will discuss the group structure of $VP_n$ and $BP_n$ in the next section.

Two virtual braid diagrams are \emph{virtually homotopic} if one diagram can be transformed into the other by a sequence of Reidemeister moves, virtual Reidemeister moves and self-crossing change. 
We will denote the set of homotopic n-strand pure virtual braids as $ H(VP_n) $ and the set of homotopic n-strand pure classical braids as $H(BP_n) $ following the notation in \cite{goldsmith}.

A \emph{representation} of an n-strand pure virtual braid is a pair 
$ (D,b) $ 
where $ D $ is a once punctured, two dimensional oriented surface with an immersed braid $b$ (where the boundary points of $b$ are contained in the boundary of $D$) modulo Reidemeister moves, Dehn twists, and handle cancellations and additions. 

\begin{figure}[htb] \epsfysize = 1 in
\centerline{\epsffile{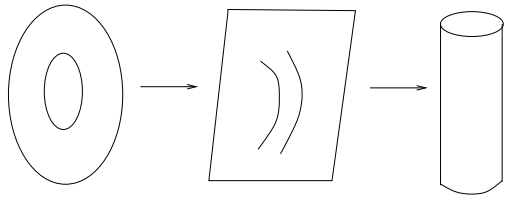}}
\caption{Surfaces for genus one representations}
\label{fig:handle}
\end{figure}

\begin{rem} We can view the surface $D$ as $ I \times I $ with $m$ attached handles ($m \geq 0 $). Theorems \ref{corres} and \ref{kup} apply to representations of pure virtual braids. In a diagram of a representation with genus one, we will draw only the handle as shown  in figure \ref{fig:handle}. Elements of $VP_n$ with genus one representations have a natural correspondence with elements of $BP_{n+2} $ and $BP_{n+1}$. 
\end{rem}

\begin{figure}[htb] \epsfysize = 1 in
\centerline{\epsffile{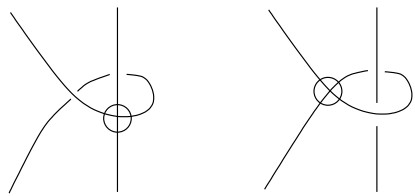}}
\caption{Braids not homotopic to a pure braid}
\label{fig:counterpure}
\end{figure}

Classical homotopy has been studied by Milnor \cite{milnor}, Goldsmith \cite{goldsmith}, and more recently Habegger and Lin \cite{xsl}.  Significant differences  exist between the classical case and the virtual case. Strikingly, not all virtual knot diagrams are homotopic to the unknot. (This topic  has been explored in \cite{d-khom} where Milnor's link groups and $ \mu $ invariants are applied to virtual link diagrams.) The difference in the case of braids is illustrated by the fact: every classical braid with a fixed ordering on the endpoints can be homotoped into a braid with no self crossings. This is not true in the virtual case, as shown in figure \ref{fig:counterpure}. In this paper, we determine which pure virtual braids are homotopic to the identity braid. 

\section{Pure  Braids}

\begin{figure}[htb] \epsfysize = 3 in
\centerline{\epsffile{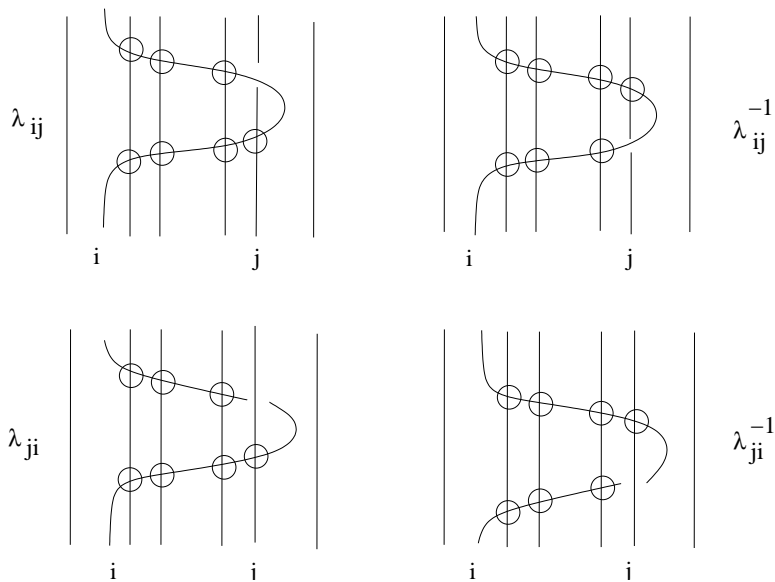}}
\caption{Generators of the pure virtual braid group}
\label{fig:pure}
\end{figure}

The structure of the n-strand pure virtual braid group is described in \cite{bardakov}. $VP_n$ is generated by the set of elements: $ \lbrace \lambda_{ik}, \lambda_{ki} | i,k \in \lbrace 1,2,...n \rbrace \rbrace $. These generators are illustrated in figure \ref{fig:pure}.

The relations in $VP_n$ are given in equations \ref{rel1} and \ref{rel2} \cite{bardakov}.
\begin{equation} \label{rel1}
 \lambda_{jk} \lambda_{in} = \lambda_{in} \lambda_{jk} \text{ for distinct } i,j,k,n \\
\end{equation}
\noindent Let $ s(ij)=1 $ if $ i<j $ and $-1 $ otherwise. Then:
\begin{equation} \label{rel2}
 \lambda_{ki} ^{s(ki)} \lambda_{kj}^{s(kj)} \lambda_{ij}^{s(ij)} = \lambda_{ij}^{s(ij)} \lambda_{kj}^{s(ij)} \lambda_{ki}^{s(ki)} 
\end{equation}
\begin{rem} The naming convention for the generators differs slightly from those given in \cite{bardakov}. In this paper, the two indices indicate the involved strands. The classical crossing is in the upper tier of crossings. The first index number indicates the over crossing strand while the second indicates the under crossing strand.
\end{rem}

\begin{figure}[htb] \epsfysize = 1 in
\centerline{\epsffile{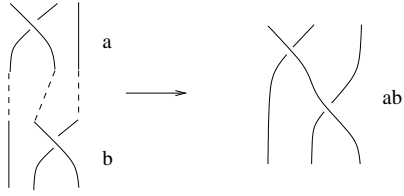}}
\caption{Braid multiplication}
\label{fig:braidmult}
\end{figure}

The n-strand pure virtual braids have a group structure where multiplication is performed by concatenating braids as shown in figure \ref{fig:braidmult}.
The n-strand pure classical braids, denoted $BP_n$, form a subgroup of $VP_n$. 
The generators of $ BP_n $ are denoted as $ \sigma_{ij} $ with $1 \leq i <j \leq n$. Each $ \sigma_{ij} $ can be expressed as a product of the virtual generators:
\begin{gather*}
\sigma_{ij} = (\lambda_{i i+1} \lambda_{i i+2}  \ldots \lambda_{i j-1}) ( \lambda_{ij} \lambda_{ji} ^{-1} ) ( \lambda_{i j-1} ^{-1} \ldots \lambda_{i i+2} ^{-1} \lambda_{i i-1} ^{-1}).
\end{gather*}
Recall the \emph{commutator} of two elements:  
\begin{equation*}
[x,y] = xy x^{-1} y^{-1}.
\end{equation*} 
Let $F_i$ denote the subgroup of $BP_n$ generated by $ \lbrace \sigma_{i i+1}, \sigma_{i i+2} \ldots \sigma_{in} \rbrace $. The set of n-strand pure classical braids homotopic to the identity is the smallest normal subgroup generated by  the commutators: 
\begin{equation}\label{classicid} 
 [\sigma_{ij} ,g \sigma_{ij} g^{-1} ] \text{ where } 
g \in F_i.
\end{equation}

\begin{figure}[htb] \epsfysize = 1 in
\centerline{\epsffile{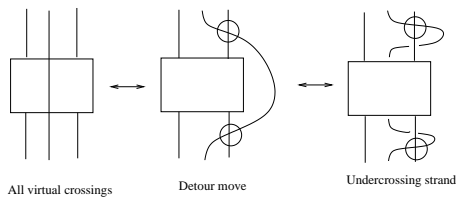}}
\caption{Converting to a braid with minimal homotopic length}
\label{fig:minlenexa}
\end{figure}

From \cite{bardakov}, we can describe $VP_n$ as a semi--direct product and give a normal form for pure virtual braids. Let $ V_n $ denote the set of generators $ \lbrace \lambda_{in}^{ \pm 1}  , \lambda_{ni}^{ \pm 1} | i \in 1,2, \ldots n-1  \rbrace $.
Then $ V_n^* $ denotes the smallest normal subgroup generated by $ V_n $ in $ VP(n-1) $. The subgroup $ V_n^* $ is normal in $ VP_n $ and  $ VP_n $ is the semi direct product:  $ V_n^* 	\rtimes  VP_{n-1} $. That is, if $ w_n \in V_n^*$ then 
\begin{equation*} 
w_n = \underset{j=1}{\overset{k}{ \prod}} g_j a_j g_j ^{-1} \text{ where } a_j \in V_n \text{ , } g_j \in VP_{n-1}.
\end{equation*} 
Based on the normal form,  we define 
the \emph{length} of $w_n$ to be $k$.
 The braid $w_n$ has \emph{minimal homotopic length} if every strand contains at least one real crossing and the braid does not contain a subsequence such that its deletion produces a braid which homotopic to $w_n$ but not virtually equivalent.  
\begin{rem}
If a braid has a strand with only virtual crossings then this braid is equivalent to a conjugate of a braid with only classical crossings on this strand.  However,  the minimal length of the braid (based on the normal form) increases as shown in figure \ref{fig:minlenexa}.
\end{rem} 

The following theorem describes a normal form for elements of $VP_n$.
\begin{thm}
Let $ b $ be an element of $VP_n $ then $ b= w_2 w_3 \ldots w_n $ where $ w_j \in V_j ^{*} $.
\end{thm}
\textbf{Proof:} See \cite{bardakov}. \qed

\begin{rem} There is a reduced form for $g_j$ when describing elements of $ V_n ^{*} $ based on the relations given earlier. \end{rem}

This result parallels Artin's theorem \cite{artin} about the normal form of pure classical braids. In the next section, we describe the set of pure virtual braids that are homotopic to the identity braid. We incorporate Goldsmith's methodology from the classical case \cite{goldsmith} and representations of virtual braids \cite{dk}.

\section{Genus and Homotopy}
We prove  that the set of n-strand pure braids homotopic 
to the identity braid is a normal subgroup of $VP_n$. We will denote this subgroup as $I(VP_n) $. 
To describe this subgroup, we will first prove a sequence of lemmas about minimal genus and $V_n ^*$. We then apply this result to a braid in normal form in the next section.

\begin{figure}[htb] \epsfysize = 3 in
\centerline{\epsffile{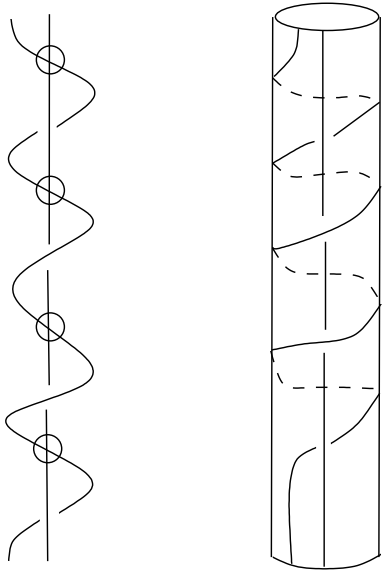}}
\caption{A representation of an element of $VP_2$}
\label{fig:v2}
\end{figure}

\begin{lem}Let $ w  \in   VP_2 $. Then a minimal genus representation of $w$ has genus less than or equal to one.
\end{lem}
\textbf{Proof:} The braid $w$ contains two strands. Immerse the second strand in a surface with one handle so that the second strand follows the longitude of the handle as shown in figure \ref{fig:v2}. \qed

\begin{cor} The elements of $VP_2 $ can  be viewed as elements of $BP_3$ and there are non-classical elements of $VP_2$ homotopic to the identity braid. \end{cor}
\textbf{Proof:} If $w$ is an element of $ VP_2$ then $w$ has a genus one representation. This representation can be viewed as an element of $BP_3$. There are non-trivial braids in $BP_3 $ that are homotopic to the identity braid. These braids in $ BP_3 $ correspond to representations of non-trivial elements of $VP_2$ which are then homotopic to the identity braid. \qed

We apply the following lemma in the cases where $n \geq 3$.
\begin{lem}\label{handle} Let $h$ be a handle in a minimal genus representation of a pure virtual braid. Let $ \lbrace s_1, s_2, \ldots s_k \rbrace  $ be the collection of strands such that the removal of $s_i$ admits a cancellation curve. Then $k \leq 2 $. \end{lem}
\textbf{Proof:}
Let $ [m] $ denote the meridian and $[l] $ the longitude of handle $h$.
Suppose that removing strand $ s_1 $ from the representation admits cancellation curve $ \alpha_1 $ on handle $h$. Then the strand $s_1 $ intersects the cancellation curve $ \alpha_1 $ transversely. By hypothesis, removing strand $ s_2 $ admits a cancellation curve $ \alpha_2 $ on $ h $. Note that $ s_2 $ intersects $ \alpha_2 $ transversely. 

Suppose that $ \alpha_2 $ is homotopic to $ \alpha_1 $ then $ s_1 $ also intersects both $ \alpha_1 $ and $ \alpha_2 $ transversely. As a result,  $ \alpha_2 $ can not be a cancellation curve after the removal of $ s_2 $. 

We observe that $ \alpha_1 $ and $ \alpha_2 $ are not homotopic but both pass through handle $h$. Assume with out loss of generality that $ \alpha_1 = [m] $ and $ \alpha_2 = [l] $ (implying that $s_1$ parallels the longitude and $s_2 $ parallels the meridian in the handle).
Let $s_3$ be a third strand in the braid whose removal admits a canceling curve $ \alpha_3 $ for $h$.  The curve $ \alpha_3 $ is homotopic to neither curve. But if $ \alpha_3 $ passes through $h$ then 
$ \alpha_3 $ is homotopic to a curve that wraps around the meridian $a$ times and the longitude $b$ times. Hence, $ \alpha_3 $ intersects at least one of the pair, $ \alpha_1 $ or $ \alpha_2 $. But then at least one the pair, $s_1 $ or $s_2$ intersects $ \alpha_3 $ transversely. Hence, $ \alpha_3 $ can not be a cancellation curve.\qed

Lemma \ref{handle} can be applied to show that certain elements of $V_n ^* $ with $ n \geq 4 $ have representations with minimal genus less than or equal to one.
 
\begin{lem}Let $n \geq 4$. If $w$ is a non-classical element of $ V_n ^* $  with minimal homotopic length then the minimal genus of a representation of $w$ is one. \end{lem}
\textbf{Proof:} Let $w$ be an element of $ V_n ^{*} $ that is homotopic to the identity such that every strand contains at least one real crossing. Suppose that a minimal genus surface for $w$ has $m$ handles. If strand $s_n$ is removed then each handle admits a canceling curve since removing $s_n$ results in the identity braid. If strand $s_i $ is removed then $m-1 $ handles admit a canceling curve. 
As a result, if some handle admits only two canceling curves then all other curves admit three or more canceling curves since $n \geq 4$. Now, by Lemma \ref{handle} there is at most one handle and the representation has genus less than or equal to one.\qed

\begin{figure}[htb] \epsfysize = 0.75 in
\centerline{\epsffile{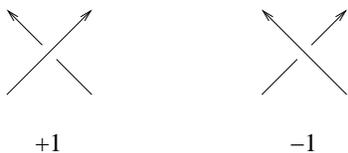}}
\caption{Sign of a classical crossing}
\label{fig:sgn}
\end{figure}

The remaining case occurs when $n=3$. We will need the following lemma about linking number. Recall that the \emph{sign} of a classical crossing, $c$, is determined by its relative orientation as shown in figure \ref{fig:sgn}.
Let $ V $ be the collection of crossings where strand $a$ passes over
strand $b$ then
\begin{equation*}
Link(a,b) = \underset{c \in V}{ \sum} sgn(c).
\end{equation*}
\begin{prop}\label{link} Let braid $b \in VP_n $ be homotopic to the n-strand identity braid. If $b$ contains $n$ copies of the generator, $ \lambda_{ij} $, then $b$ contains $n$ copies of generator:  $\lambda_{ij} ^{-1} $. \end{prop}
\textbf{Proof: }For the generator $ \lambda_{ij} $, $Link(i,j) = -1 $ and $ Link(j,i) = 1 $. Linking number is a homotopy invariant and in the identity braid $ Link(i,j) = 0 $.  Hence the generator $ \lambda_{ij}$ paired with the generator $ \lambda_{ij} ^{-1}.$\qed

Let $ \lambda_{p(ij)} $ represent either $ \lambda_{ij} $ or $ \lambda_{ji} $. We use this notation in the following lemma.

\begin{lem} If $w$ is a braid in $ V_3 ^ * $ with minimal homotopic length that is homotopic to the identity then the minimal genus of a representation of $ w$ is less than or equal to one. \end{lem}

\textbf{Proof:}
Let $w$ be a braid in $V_3 ^*$ with minimal homotopic length that is homotopic to the identity braid.  
Suppose that a minimal genus representation of $w$ is a surface with $m$ handles. 
Consider the representation of $w$.
Since $w$ is in normal form, the removal of the $3^{rd}$ strand from the representation results in the 2-strand identity braid and $m$ canceling curves on the surface. 

The removal of the first strand results in the braid $x$, a 2-strand braid that is homotopic to the identity and as a result, $m-1$ canceling curves on the surface. 
Similarly, the removal of the $2^{nd}$ strand results in the braid $y$, a 2-strand braid homotopic to the identity and $m-1 $ canceling curves on the surface. 

Applying Lemma \ref{handle}, we observe that the representation of $w$ contained at most 2 handles. 
Note that if $x$ is classical, then $x$ is the identity braid and the representation of $w$ has genus less than or equal to one.  
As a result, we will assume that both $x$ and $y$ are non-classical for the remainder of the proof. 

Let $\hat{x}$ denote the subsequence of $w$ consisting of conjugates of $ \lambda_{p(13)} ^{\pm 1}$. The removal of the second strand from $ \hat{x}$ results in a 2-strand braid homotopic to the identity, while the removal of the $1^{st} $ and $3^{rd}$ strands results in the identity braid. 
Hence, a minimal genus representation of $ \hat{x}$ has genus one. 
(Note that if the genus is zero, then $ \hat{x}$ is a classical braid which contradicts our assumption that $ x $ is a non-classical braid.)
Let $ \hat{y} $ denote of the subsequence of $w$ consisting of conjugates of $ \lambda_{p(23)} ^{\pm 1}$. Following the argument given for $\hat{x} $, a minimal genus representation of $ \hat{y} $ also has genus one.

\begin{figure}[htb] \epsfysize = 2 in
\centerline{\epsffile{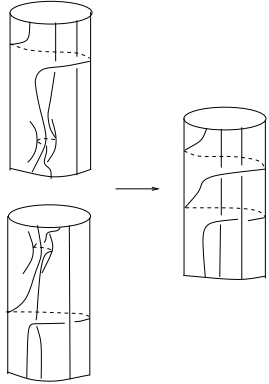}}
\caption{Minimal genus surface}
\label{fig:min3}
\end{figure}

Now, $w= x_1 y_1 x_2 y_2 \ldots x_n y_n $ where $\hat{x}= x_1 x_2 \ldots x_n $ and $ \hat{y} = y_1 y_2 \ldots y_n $. Because minimal genus representations of $\hat{x} $ and $\hat{y}$ have genus one then we may assume that (after isotopy) minimal genus representations of each subbraid $x_i$ and $y_i$ occur on a tube (possibly) with handles as shown in figure \ref{fig:min3}.

 However, a representation of $w$ has at most genus two. If handles occur in a representation of some $x_i$ (or $y_i$) then all the handles are canceled in a representation of $\hat{x}$ (or $\hat{y}$). Without loss of generality, assume that the tube representation of some $ x_i$ contains at least one handle. (The handles can be selected so that a handle either contains the strands $1$ and $3$ or the strands $1$ and $2$.) Suppose that this handle contains the strands $1$ and $3$ then
the subbraid $ x_{i+1}$ must have a corresponding handle involving these strands. 
 
In the braid $w$, the braid $y_i$ occurs between $x_i$ and $x_{i+1}$.  As a result, the handle in $x_i$ cancels with a handle occurring in either $y_{i-1}$ or $y_i$. But $ \hat{y} $ contains no virtual or classical crossings between $1 $ and $3$. Hence, the handles can not be canceled in $w$ and the minimal genus of $w$ is three, a contradiction.

Suppose that the handles in $x_i$ that can not be removed in $w$ contain the strands $1$ and $2$. That is, $g_i^{-1} $ does not cancel with $g_{i+1}$. Now, the deletion of either strand one or strand two admits a canceling curve for each of these handles, but not the tube. The deletion of strand three admits a canceling curve for every handle. Again by applying Lemma \ref{handle}, we observe that the representation of $w$ has genus one.\qed

We have proved the following theorem:
\begin{thm} \label{min-g} Let $w_n$ be a minimal length, non-trivial element of $V_n ^*$ such that every strand contains a real crossing. If $w_n$ is homotopic to the identity braid then a minimal genus representation of $w_n$ has genus less than or equal to one.\end{thm}

\begin{rem}If a braid contains a strand with only virtual crossings then some representation contains a handle with this strand immersed along the longitude.  \end{rem}
	
\begin{lem} Let $ b= w_2 w_3 \ldots w_n $ be an element of $ VP_n $ (where $ w_j \in V_j ^* $). Then $b $ is homotopic to the n-strand identity braid if and only if each $ w_i $ is homotopic to the identity braid. \end{lem}

\textbf{Proof: }  If each $w_i $ is an element of $I(VP_n)$ then 
$ b = w_1 w_2 \ldots w_n $ is an element of $I(VP_n) $. Let $w$ be an element of $I(VP_n) $. Then there is a homotopy sequence $ \lbrace p_0, p_1, \ldots p_m \rbrace $, with $ w = p_0 $ and $p_m$ equivalent to the identity braid, that transforms $w $ into the identity. Remove strand $n$ from each diagram and replace it with the identity strand. This reduces $w$ to the braid $b_{n-1} $ where $b_{n-1}= w_1 w_2 \ldots w_{n-1} $. This preserves the homotopy sequence so that $b_{n-1} $ is homotopic to the identity strand. Removing strands, we observe that $b_i = w_1 w_2 \ldots w_i $ is homotopic to the identity braid. Hence, each $w_i$ is homotopic to the identity braid.\qed

\section{Braids homotopic to the identity braid}
\begin{figure}[htb] \epsfysize = 5 in
\centerline{\epsffile{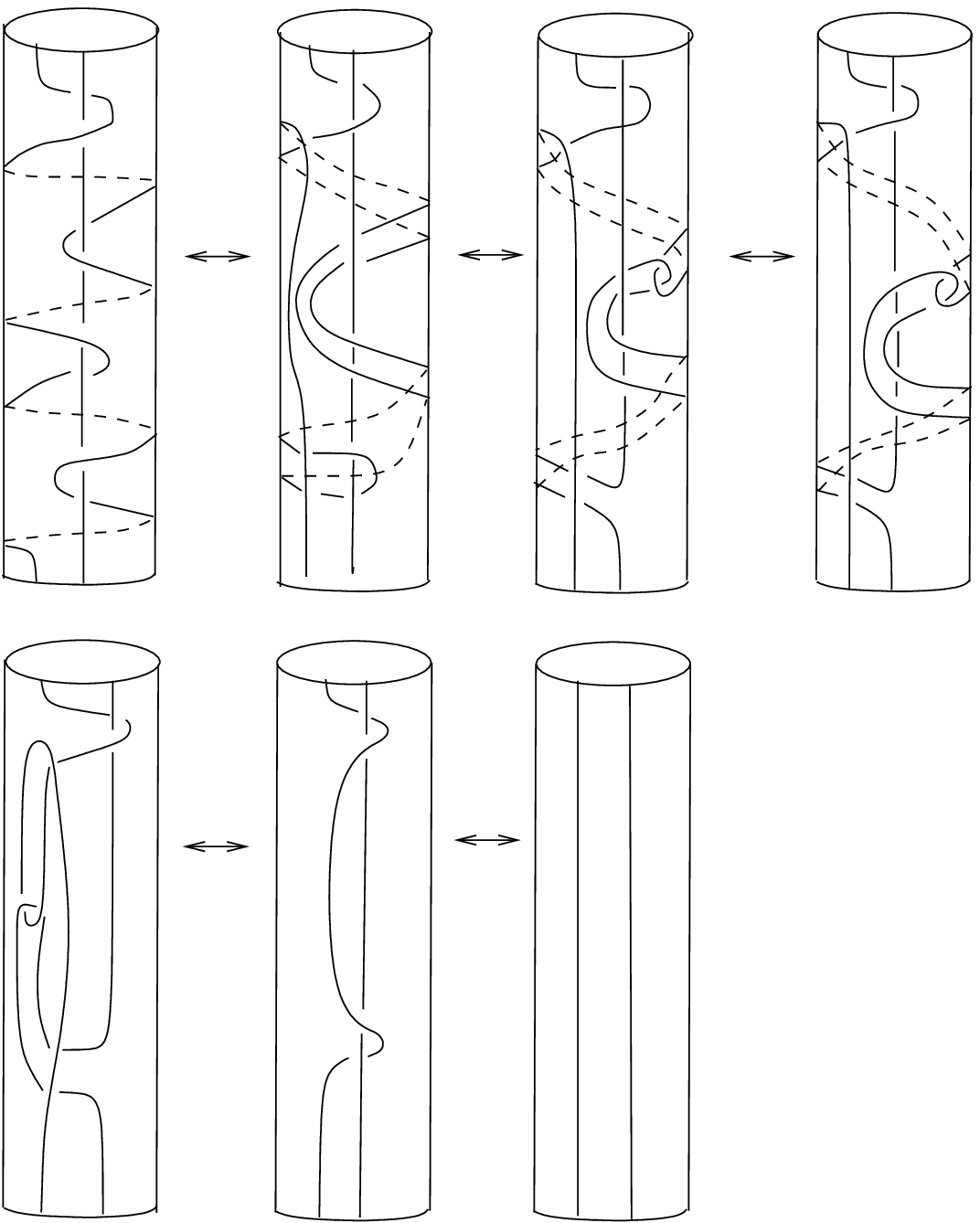}}
\caption{A homotopy sequence in $VP_2$}
\label{fig:VP2homotopy}
\end{figure}

\begin{figure}[htb] \epsfysize = 5 in
\centerline{\epsffile{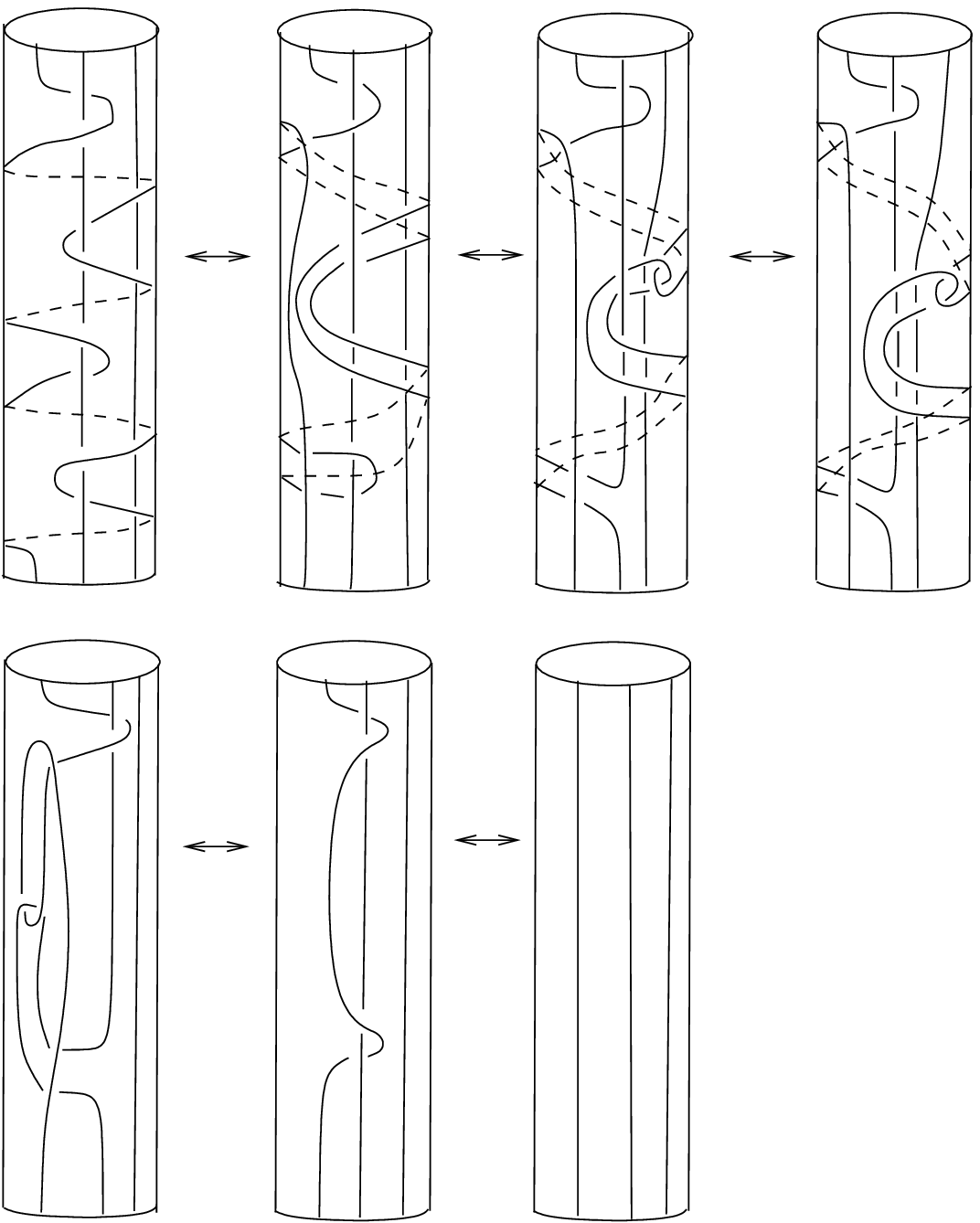}}
\caption{A homotopy sequence in $VP_3$}
\label{fig:VP3homotopy}
\end{figure}

\begin{figure}[htb] \epsfysize = 5 in
\centerline{\epsffile{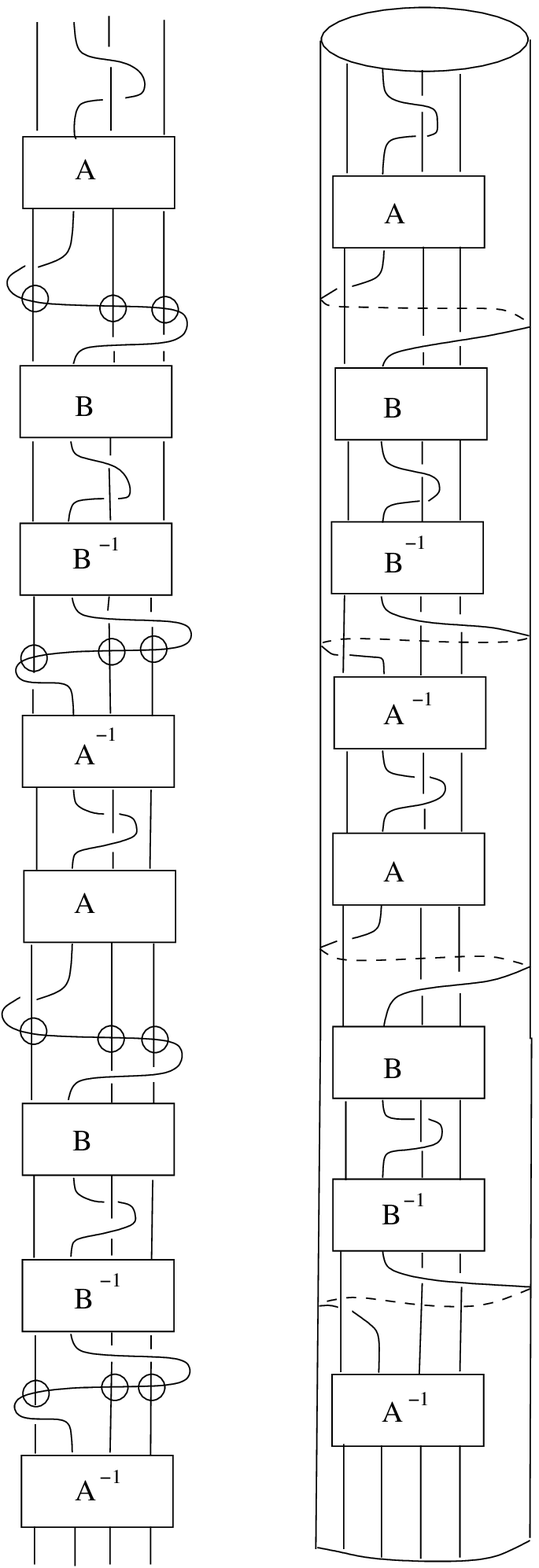}}
\caption{A representation of $w$}
\label{fig:brep}
\end{figure}

We show that the set of braids homotopic to the identity braid form a normal subgroup of $VP_n$. Let $x$ be an element of $VP_n$ such that $ x = (\lambda_{i-1 i} \lambda_{i-2 i} \ldots \lambda_{1 i}) (\lambda_{ni}^{-1} )(\lambda_{n-1 i} ^{-1} \ldots \lambda_{i+1 i} ^{-1} )$.  Let $g_a $ and $g_b $ denote classical braids generated by the set: $ \lbrace \sigma_{1i}, \sigma_{2i} \ldots \sigma_{i-1 i} \sigma_{ii+1} \ldots \sigma_{in} \rbrace $.
\begin{prop}If an element of $VP_n$ has the form $ [\sigma_{ij}, g_a x g_b \sigma_{ij} g_b ^{-1} x^{-1} g_a ^{-1} ] $  (as shown in figure \ref{fig:brep} where $g_a $  and $g_b$ are denoted as $A$ and $B$) then the braid is homotopic to the identity braid. \end{prop}
\textbf{Proof:}  
We show a sample homotopy sequences  for elements of $VP_2$ and $VP_3 $  in figures \ref{fig:VP2homotopy}. 
and \ref{fig:VP3homotopy}. \qed

\begin{rem} This braid can be expressed as the product of an element of $V_n ^* $ and $ VP_{n-1}$. \end{rem}

\begin{prop} If $w$ is an  element of $ V_n ^* $ that is homotopic to the identity with minimal homotopic length then some representation of $w$ corresponds to an element of $BP_{n+1}$ with the form:
$ [ \sigma_{ij} , g \sigma_{ij} g^{-1} ] $. \end{prop}

\textbf{Proof:} Let $w$ be a braid in $V_n ^ * $ with minimal homotopic length that is homotopic to the identity. By Theorem \ref{min-g}, the braid $w$ has a representation with minimal genus one, a torus. View the torus as the complement of two linked curves to obtain an element of $BP_{n+2} $. After removing strand $n+2 $, we obtain a classical braid in $BP_{n+1}$ that is homotopic to the identity. The braid $w$, as an element of $BP_{n+1} $,  can written in the classical normal form, $b_1 b_2 \ldots b_{n+1} $, where $b_i$ is an element of the group generated by $ \lbrace \sigma_{i i+1}, \sigma_{i i+2}, \ldots \sigma_{i n+1} \rbrace $. In $BP_{n+1}$, a braid
 homotopic to the identity has the form: $ [ \sigma_{i j} , g \sigma_{ij} g^{-1} ]$.
  
In a representation of $w$ the longitude of the torus corresponds to strand $n+1$. As a result, consider the first $b_i $ that includes the $n+1^{th}$ strand. 

\begin{figure}[htb] \epsfysize = 0.75 in
\centerline{\epsffile{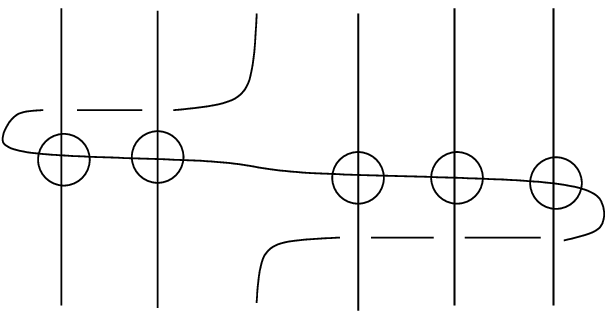}}
\caption{Braid x in $VP_n$}
\label{fig:xbraid}
\end{figure}

\begin{figure}[htb] \epsfysize = 0.75 in
\centerline{\epsffile{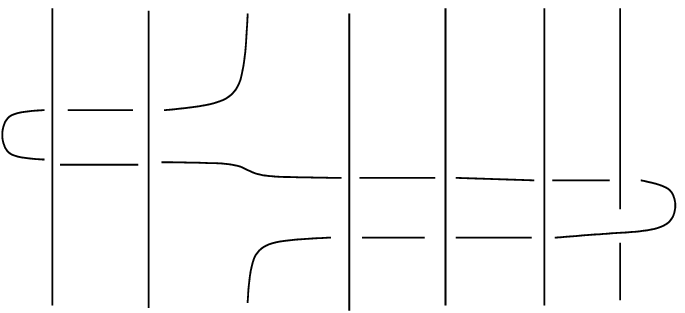}}
\caption{Braid $ \hat{x} $ in $BP_{n+1} $}
\label{fig:hatxbraid}
\end{figure}

Now, $ \sigma_{i n+1} $ either holds the place of $ \sigma_{ij} $  or is term in $g$ from equation \ref{classicid}.  Let $ x $ denote the element of  $VP_n$ shown in figure \ref{fig:xbraid}.
In $VP_n$, 
\begin{equation*}
x = (\lambda_{i-1 i} \lambda_{i-2 i} \ldots \lambda_{1 i}) (\lambda_{ni}^{-1} \lambda_{n-1 i} ^{-1} \ldots \lambda_{i+1 i} ^{-1} ).
\end{equation*}
From a representation of $x$, we obtain $\hat{x}$, the element of 
$BP_{n+1}$ shown in figure \ref{fig:hatxbraid}.
In $BP_{n+1}$, 
\begin{equation*}
\hat{x} = (\sigma_{i i+1} ^{-1} \sigma_{i i+2} ^{-1}  \ldots \sigma_{i n} ^{-1})\sigma_{i n+1} ^{-1} ( \sigma_{i n} \ldots \sigma_{i i+1} ) .
\end{equation*} 
If $ \sigma_{i n+1} $ corresponds to $ \sigma_{ij} $ then the homotopy sequence requires a crossing change on strand $n+1$. 
Hence $g$ must contain the term $ \sigma_{i n+1} $ and the braid $b_i $ has the form:
\begin{equation*}
[ \sigma_{ij} , q_a \hat{x} q_b \sigma_{ij} q_b ^{-1} \hat{x}^{-1} q_a^{-1} ] 
\end{equation*}
where $ q_a $ and $q_b $ are classical braids in the group  generated by
$ \lbrace \sigma_{i i+1} \ldots \sigma_{i n} \rbrace$.
This corresponds to the braid shown in figure \ref{fig:brep}. In terms of $VP_n$ this is a braid of the form:
\begin{equation}\label{almost}
[\sigma_{ij}, G_a \lambda_{ni} ^{-1} G_b \sigma_{ij} G_b ^{-1} \lambda_{ni} G_a^{-1} ]
\end{equation}
(Note that $G_a $ and $G_b$ are the images of $q_a $ and $q_b$ with the appropriate part of $\hat{x} $ in $VP_n$.)
Let 
\begin{align} \label{ende}
\lambda &= \lambda_{ni} \\
 z &= G_a ^{-1} \sigma_{ij} G_a  \label{ende2} \\
 y &= G_b \sigma_{ij} G_b  \label{ende3}.
 \end{align}
Then the braid in equation \ref{almost} can be rewritten as:
\begin{equation*}
G_a ( z \lambda ^{-1} y \lambda z^{-1} \lambda ^{-1} y ^{-1} \lambda ) G_a ^{-1}
\end{equation*}
This can be rewritten in normal form for $VP_{n}$:
\begin{equation*}
G_a (z y z^{-1} y^{-1} ) (yzy^{-1} \lambda ^{-1} y z^{-1} y^{-1})
(yz \lambda z^{-1} y^{-1})(y \lambda^{-1} y^{-1}) \lambda G_a ^{-1} 
\end{equation*}
The braid $ z y z^{-1} y^{-1} $ is an element of $VP_{n-1}$ and homotopic to the identity. Then
$yzy^{-1} $ is homotopic to the braid $z$. Canceling terms, we obtain the braid:
\begin{equation}
G_a z \lambda ^{-1} y \lambda z^{-1} \lambda^{-1} y^{-1} \lambda G_a ^{-1}.
\end{equation}
Using equations \ref{ende}, \ref{ende2}, and \ref{ende3} to rewrite, we obtain:
\begin{equation}
\sigma_{ij} G_a \lambda ^{-1} G_b \sigma_{ij} G_b ^{-1} \lambda G_a^{-1} \sigma_{ij} G_a \lambda ^{-1} G_b \sigma_{ij} G_b ^{-1}   
\lambda G_a ^{-1}
\end{equation}
This is the same form as the original braid.\qed 
\begin{rem} Note that in this form, the $i^{th} $ strand underpasses all $i-1$ previous strands. Any weaving on these strands can be described by multiplying $w$ by the braid $b_k $ with $k \leq i $. Note that if
the original braid $w$ over crosses the $i-1$ previous strands, then $b_k$ with $k < i$ involves crossings on the $n+1^{th}$ strand. This contradicts the fact that $b_i$ is the first strand to have crossings involving strand $n+1$. \end{rem}

We have proved the following theorem:
\begin{thm} The set of pure virtual braids homotopic to the identity is the smallest normal subgroup generated by elements of the form:
\begin{equation*}
[ \sigma_{ij}, g_a x g_b \sigma_{ij} g_b ^{-1} x ^{-1} g_a ^{-1} ] 
\end{equation*}
where $ x = (\lambda_{i-1 i} \lambda_{i-2 i} \ldots \lambda_{1 i}) (\lambda_{ni}^{-1} \lambda_{n-1 i} ^{-1} \ldots \lambda_{i+1 i} ^{-1} ) $ or $x$ is the identity braid and $ g_a$, $g_b$ are classical braids generated by the set:
$ \lbrace \sigma_{1i}, \sigma_{2i}, \ldots \sigma_{i-1 i}, \sigma_{i i+1}  \ldots \sigma_{i n} \rbrace $. \end{thm}

\noindent \textbf{Acknowledgments: }This research was performed while the author held a National Research Council Research
  Associateship Award jointly at the Army Research Laboratory and the U.\S.\ Military Academy. This effort was sponsored in part by the
Office of the Dean at the United States Military Academy.
 The U.S. Government is authorized to reproduce and distribute
reprints
for Government purposes notwithstanding any copyright annotations thereon.
The
views and conclusions contained herein are those of the authors and should
not be
interpreted as necessarily representing the official policies or
endorsements,
either expressed or implied, of the United States Military Academy or the U.S. Government. (Copyright
2008.)


\begin{thebibliography}{10}

\baselineskip=12pt % space between lines
\parskip=2pt plus 1pt % space between paragraphs

\bibitem{artin}
\newblock{Artin, E. Theory of braids. Ann. of Math. (2) 48, (1947). 101--126}

\bibitem{bardakov}
\newblock{Bardakov, Valerij G. The virtual and universal braids. Fund. Math. 184 (2004), 1--18}

\bibitem{dk}
\newblock{Dye, H. A. and Kauffman, Louis H.  Minimal surface representations of virtual knots and links. Algebr. Geom. Topol. 5 (2005), 509--535 }

\bibitem{d-khom}
\newblock{Dye, H. A. and Kauffman, Louis H. Virtual Homotopy. Preprint: www.arxiv.org/GT.}


\bibitem{goldsmith}
\newblock{Goldsmith, Deborah Louise.
Homotopy of braids---in answer to a question of E. Artin. Topology Conference (Virginia Polytech. Inst. and State Univ., Blacksburg, Va., 1973), 91--96. Lecture Notes in Math., Vol. 375, Springer, Berlin, 1974}

\bibitem{xsl}
\newblock{Habegger, Nathan and Lin, Xiao-Song. The classification of links up to link-homotopy. J. Amer. Math. Soc. 3 (1990), no. 2, 389--419}

\bibitem{kamadabraid}
\newblock{Kamada, Seiichi. Invariants of virtual braids and a remark on left stabilizations and virtual exchange moves. Kobe J. Math. 21 (2004), no. 1-2, 33--49. } 

\bibitem{kamada2}
\newblock{Kamada, Naoko and  Kamada, Seiichi. Abstract link diagrams and virtual knots. J. Knot Theory Ramifications 9 (2000), no. 1, 93--106}

\bibitem{kvirt}
\newblock{Kauffman, Louis H. Virtual Knot Theory. European Journal of Combinatorics, Vol. 20 (1999), No. 7, 663--690,}

\bibitem{sophia}
\newblock{Kauffman, Louis H. and Lambropoulou, Sofia. Virtual braids. Fund. Math. 184 (2004), 159--186.}



\bibitem{kup}
\newblock{Kuperberg, Greg. What is a virtual link? Algebr. Geom. Topol. 3 (2003), 587--591 }


\bibitem{milnor}
\newblock{Milnor, John. Link groups. Ann. of Math. (2) 59, (1954). 177--195}




\end{thebibliography}
\end{document}